 \documentclass{article}
 \usepackage{amssymb}
 \usepackage{amsmath}
 \newcommand{\be}{\begin{equation}}
 \newcommand{\ee}{\end{equation}}
 \newcommand{\bd}{\begin{displaymath}}
 \newcommand{\ed}{\end{displaymath}}
 \newcommand{\bea}{\begin{eqnarray}}
 \newcommand{\eea}{\end{eqnarray}}
 \newcommand{\beas}{\begin{eqnarray*}}
 \newcommand{\eeas}{\end{eqnarray*}}
 \newcommand{\bc}{\begin{center}}
 \newcommand{\ec}{\end{center}}
 \newcommand{\rr}{\mathbb{R}}

 \input{amssym.def}
 \input{amssym.tex}

\begin{document}
 \title{ STABILITY OF SOLITARY WAVES FOR THREE COUPLED
             LONG WAVE - SHORT WAVE INTERACTION EQUATIONS}

 \author
 { Handan Borluk\footnote{e-mail: hborluk@isikun.edu.tr} ~and~
 Saadet Erbay\footnote{e-mail: serbay@isikun.edu.tr} \\
 Department of Mathematics,
 Isik University, 34980 Sile-Istanbul, Turkey}

 \maketitle

 \begin{abstract}
 \noindent
In this paper we consider a three-component system of one dimensional long wave-short wave
interaction equations. The system has two-parameter  family of solitary wave solutions.
We prove  orbital stability of  the solitary wave solutions using variational methods.
 \end{abstract}



 \newpage

 \section{Introduction}

 \noindent
 In the present paper, we study orbital stability  of solitary wave
 solutions  for the three coupled  long wave-short wave
 interaction (LSI) equations
  \be
 \left . \begin{array}{ll}
  & i\phi_t+ \phi_{xx}=\beta w\phi,  \\
  &  i\psi_t+ \psi_{xx}=\beta w\psi, \\
  & w_t= \beta (|\phi|^2+|\psi|^2)_x,
 \end{array} \right \} \label{lwsw}
 \ee
 where
 $\phi,\psi: \rr \times {\rr^+} \rightarrow {\mathbb{C}},
 ~w: \rr \times {\rr^+} \rightarrow \rr $,
 $\beta$ is a real constant.  Here $w$ represents  a long wave mode, and
 $\phi$ and $\psi$ denote short wave modes propagating in a continuous  medium.
 This system describes the resonant interaction among two short wave modes  with equal
 group speeds and one long wave mode whose phase speed is equal to
 the group speed of short waves. System (\ref{lwsw}) appears, for instance,  in water
 waves \cite{ma, craik} and in a bulk elastic medium \cite{erbay}. In a recent study
 \cite{borluk1}, it is shown that system (\ref{lwsw}) has a two-parameter family of solitary
 wave solutions of the form
 \be
 \left . \begin{array}{ll}
 & \phi_s(x,t)=\Phi(x-ct)e^{i \omega t}, \\
 & \psi_s(x,t)=\Psi(x-ct)e^{i \omega t}, \\
 & w_s(x,t)= W(x-ct),
 \end{array} \right \} \label{solution}
 \ee
 where  $\displaystyle W(x)=-\beta (|\Phi(x)|^2+|\Psi(x)|^2)/c$,
 $~\displaystyle (\Phi(x),\Psi(x))=(R_1(x), R_2(x))e^\frac{i cx}{2}$,
 $c>0$ and $4\omega -c^2>0$. Here $W\in L^2(\rr)$ and $(R_1, R_2)\in H^1(\rr)\times H^1(\rr) $
 are positive solutions of
 \be
  \left . \begin{array}{ll}
 & \displaystyle -u_{xx} +(\omega-\frac{c^2}{4}) u- \frac{\beta^2}{c}  (u^2+v^2)u=0,  \\
 & \displaystyle -v_{xx} +(\omega-\frac{c^2}{4}) v- \frac{\beta^2}{c}  (u^2+v^2)v=0.
 \end{array} \right \} \label{cs1}
 \ee
 System (\ref{lwsw}) is a generalization of the two component  long wave-short wave
 interaction system given by
 \be
 \left . \begin{array}{ll}
 &  i\phi_t+ \phi_{xx}=\beta w\phi,  \\
 & w_t= \nu (|\phi|^2)_x.
 \end{array} \right \} \label{lwsw2}
 \ee
 System (\ref{lwsw2}) was derived  to describe  the resonant interaction of a short
 wave, $\phi$, and a long wave, $w$, propagating on the surface of water  \cite{red}.
 The same system was also obtained  for the resonant interaction of internal
 gravity waves \cite{grim}.  In \cite{laurencot}, orbital stability of solitary
 wave solutions of the two component  system  (\ref{lwsw2}) was considered and
 it was observed that the solitary waves of the form
 \bd
  \phi(x,t)=R(x-ct) ~e^{ i\omega t+ i\frac{c}{2}(x-ct) }, \;\;\;
  w(x,t)=W(x-ct),
 \ed
 where $(R, W) \in H^1(\rr)\times L^2(\rr)$,
 are orbitally stable when  $c>0$ and $4\omega -c^2>0$.  In that study, a Lyapunov
 functional, whose critical point is the solitary wave solution of (\ref{lwsw2}),
 was constructed using the conserved quantities of the system.  It was shown that
 the stability of solitary waves relied on suitable lower and upper bounds on the
 variation of the Lyapunov functional. This method has been developed in  \cite{ben}
 to prove the stability of solitary waves of the Korteweg -de Vries equation.
 In a later study \cite{wein1}, the same method has been used to show the stability of
 standing waves of the Nonlinear Schr\"odinger equation, which has been already proved in
 \cite{caze} using the concentration-compactness methods. In the present paper, we prove
 that solitary waves (\ref{solution})   are  orbitally stable   for the LSI system
 (\ref{lwsw}) when $c>0$ and $4\omega -c^2>0$, using a variational method, the so-called
  Lyapunov method.

 The organization of the paper is as follows: The local well-posedness of  the Cauchy problem
 for (\ref{lwsw}) is discussed, and conserved integrals for the  same system is given in
 section 2. A variational characterization of the solitary waves,
 which will be used  in the proof of the stability of solitary wave solutions, is briefly
 presented in section 3. We state the stability theorem
 that relies on a lower  bound of the second variation of the Lyapunov
 functional in section 4. Using the analysis of the unconstrained variational
 problem carried out in \cite{borluk1}, the lower bound is proved and
 the stability of solitary waves is established in the same section.

 \emph{Notations.} Throughout the paper $L^p(\rr),~1\leq p <\infty$, represents the space of
 $p-$integrable functions. $\|f\|_p$ denotes the $L^p(\rr)$ norm of $f,~1\leq p\leq \infty$.
 $H^1(\rr)=W^{1,2}(\rr)$ is the Sobolev space of $f$  for which the norm
 $\|f\|_{H^1}^2=\|f\|_2^2+\|\nabla f\|_2^2$ is finite. $\langle f, g
 \rangle$ refers to the inner product of $f$ and $g$ in $L^2(\rr)$.

 \setcounter{equation}{0}
 \section{Local Well-Posedness of Cauchy Problem}

  \noindent The Cauchy problem for the two component LSI system (\ref{lwsw2}) was studied in  \cite{tsutsumi} for initial
 data $(\phi_0,w_0) \in H^{1/2}(\rr)\times L^2(\rr)$. A contraction technique  was used to prove existence and
 uniqueness of solutions  of the initial value problem in suitable Banach spaces. Conservation of energy and
 masses were used to extend the local solution globally. Later, the local-well-posedness result for (\ref{lwsw2})
 was improved in \cite{ginibre} for initial data $(\phi_0,w_0) \in H^k(\rr)\times L^{1/k}(\rr),~~0<k<1/2$.
 The Cauchy problem for the three component LSI system (\ref{lwsw}) was considered in  \cite{borluk3} for
 initial  data $(\phi_0,\psi_0, w_0) \in H^{1/2}(\rr)\times H^{1/2}(\rr)\times L^2(\rr)$. Following the ideas of
 \cite{tsutsumi}, the local-well-posedness theory was established in the same study. In the present section,
 the local well-posedness of the Cauchy problem for (\ref{lwsw}), which is necessary in the study of stability
 of solitary waves, will be presented briefly.

 For the three component LSI system, the Cauchy problem
 \be
 \left . \begin{array}{ll}
 & i\phi_t+ \phi_{xx}=F_\phi(\phi,\psi)  \\
 & i\psi_t+ \psi_{xx}=F_\psi(\phi,\psi) \\
 & \phi(x,0)=\phi_0(x),~~ \psi(x,0)=\psi_0(x)
 \end{array} \right \} \label{lwsw21}
 \ee
 was considered in \cite{borluk3},  where
 \bd
 (F_\phi ,F_\psi)= \left (\beta^2\int\limits_0^t
        (|\phi(x,s)|^2+|\psi(x,s)|^2)_x d s +\beta w_0(x) \right ) ~(\phi(x,t),\psi(x,t)).
 \ed
 Here $(\phi_0, \psi_0, w_0) \in H^{1/2}(\rr)\times H^{1/2}(\rr)\times L^2(\rr)$.
 Following \cite{tsutsumi}, a fixed point method was used to prove existence and
 uniqueness  of local in time solutions. First, the Cauchy problem (\ref{lwsw21}) was
 written as a coupled system of integral equations
 \be
 \left . \begin{array}{ll}
 &  \phi(t)=J_1(\phi,\psi)=U(t) \phi_0(x)-i \int\limits_0^t  U(t-s)F_\phi(\phi(s), \psi(s)) d s\\
 &  \psi(t)=J_2(\phi,\psi)=U(t) \psi_0(x)-i \int\limits_0^t  U(t-s)F_\psi(\phi(s), \psi(s)) d s,
 \end{array} \right \} \label{lwsw22}
 \ee
 where $U(t)=e^{it \frac{\partial^2}{\partial x^2}}$ is the Schr\"odinger linear group.
 Similar to that of \cite{tsutsumi}, the function space $X^T$ was
 defined by
 \bd
 X^T=\{ f: [0,T]\times \rr \rightarrow \mathrm{C}| f\in C([0,T]; H^{1/2}(\rr)),~
        f_x\in L^\infty (\rr; L^2[0,T])  \},
 \ed
 and endowed with the norm
 \bd
 \|f\|_{X^T}=\sup\limits_{0\leq t\leq T}\|f(.,t)\|_{H^{1/2}}
            +\sup\limits_{\rr} \left (\int\limits_0^T |f_x(x,t)|^2 dt \right )^{1/2}.
 \ed
 Thus $(\phi,\psi)\in Y^T=X^T\times X^T$ for which the norm was defined as
 \be
 \|(\phi,\psi )\|_{Y^T}=\|\phi\|_{X^T}+ \|\psi\|_{X^T}.
 \ee
 A ball of radius $R>0$  was defined as
  $ B_{R}= \{ (\phi,\psi)\in {Y^{T}}: \| \phi \|_{X^T}+\| \psi \|_{X^T} \leq R \}$.
  Using smoothing effects estimates (Lemmas 1-3 in \cite{tsutsumi}) obtained in \cite{kenig1, kenig2}
  and  Lemmas 4-5 in \cite{tsutsumi}, the norm of  $\Psi(\phi,\psi)= (J_1, J_2)$
 in $Y^T$ was calculated to show $\Psi:  B_R \rightarrow B_R$:
 \beas
 \| \Psi(\phi,\psi) \|_{Y^T} &=&  \|J_1(\phi,\psi)\|_{X^T}+\|J_2 (\phi,\psi)\|_{X^T}\\
 &  \leq  &   C_1 \|\phi_0\|_{H^{1/2}}+ C_2 \|\psi_0\|_{H^{1/2}}  
    + K_1(T) (\|\phi\|_{X^T}^2\\
    &&+\|\psi\|_{X^T}^2+1) \|\phi\|_{X^T} 
    + K_2(T) (\|\phi\|_{X^T}^2+\|\psi\|_{X^T}^2+1) \|\psi\|_{X^T}\\
 & \leq & C( \|\phi_0\|_{H^{1/2}}+ \|\psi_0\|_{H^{1/2}}) \\
&& + K(T) (\|\phi\|_{X^T}^2+\|\psi\|_{X^T}^2+1)
 (\|\phi\|_{X^T}+ \|\psi\|_{X^T}).
 \eeas
 Here $C=\max\{C_1,C_2\}$ is a positive constant,  $K(T)=\max\{K_1(T), K_2(T)\}~$ and
 $ K(T)\rightarrow 0$ as  $T\rightarrow 0$. $R$ was chosen so as to
 $ ~C (\|\phi_0\|_{H^{1/2}}+\|\psi_0\|_{H^{1/2}}) \leq  R/2$. Then for fixed $R>0$, small $T$ was
 taken to ensure $~K(T)(4R^3+2R)< R/2$. This led to
\bd
 \| \Psi(\phi,\psi) \|_{Y^T} \leq R.
\ed

 For uniqueness of the fixed point, in a similar way, the norm of
 $\Psi(\phi_1, \psi_1)(t)-\Psi(\phi_2, \psi_2)(t)$ in $Y^T$  was calculated:
 \beas
  &&\| \Psi(\phi_1, \psi_1)(t)-\Psi(\phi_2, \psi_2)(t) \|_{Y^T}=
    \| J_1(\phi_1,\psi_1)- J_1(\phi_2,\psi_2) \|_{X^T}\\
 &&~~~~~~~~~~~~~~~~+ \| J_2(\phi_1,\psi_1)- J_2(\phi_2,\psi_2)\|_{X^T}\\
 && ~~~~~~~~~~~~~~~~\leq H_1(T) \left ( \|\phi_1\|^2_{X^T}+ \|\phi_2\|^2_{X^T}+ \|\psi_1\|^2_{X^T}\right.\\
 && ~~~~~~~~~~~~~~~~ \left.+\|\psi_2\|^2_{X^T} +1 \right ) \|\phi_1-\phi_2\|_{X^T} \\
 && ~~~~~~~~~~~~~~~~+ H_2(T) \left ( \|\psi_1\|^2_{X^T}+ \|\psi_2\|^2_{X^T}+ \|\phi_1\|^2_{X^T}\right.\\
 && ~~~~~~~~~~~~~~~~+ \left.\|\phi_2\|^2_{X^T} + 1\right ) \|\psi_1-\psi_2\|_{X^T}\\
 && ~~~~~~~~~~~~~~~~\leq H(T) (4R^2+1) \|(\phi_1-\phi_2,\psi_1-\psi_2 )\|_{Y^T},
 \eeas
 where  $H(T)=\max\{H_1(T), H_2(T)\}$ and  $H(T) \rightarrow 0$ as  $T\rightarrow 0$.
 $T$ was chosen as a small quantity so that  $ H(T)(4R^2+1)<1$. It was concluded in \cite{borluk3}  that
 $\Psi(\phi,\psi)$ was a contraction
 on $B_R$ and $(\phi, \psi)$ was the unique solution of (\ref{lwsw22}), i.e. the local well-posedness of the
 Cauchy problem for (\ref{lwsw}) was established:

\textbf{Theorem}
Let $(\phi_0,\psi_0) \in  H^{1/2}(\rr)\times H^{1/2}(\rr)$ and $w_0
\in L^2(\rr)\cap L^\infty(\rr)$. There exists a unique solution
$\left(\phi(x,t),\psi(x,t)\right)$ of the Cauchy problem
(\ref{lwsw21}) on $[0,T]$ for $T>0$ such that $\phi\in C\left([0,T];H^{1/2}(\rr)\right)$,
$\phi_x\in L^\infty \left(\rr;L^2[0,T]\right)$,
$\psi\in C\left([0,T];H^{1/2}(\rr)\right)$ and
$\psi_x\in L^\infty \left(\rr; L^2[0,T]\right)$.

 The conserved integrals of the LSI system (\ref{lwsw}) are of the form \cite{borluk2}
 \bea
 && I_1=\int\limits_{\rr}|\phi|^2~d x,~~~~
    I_2=\int\limits_{\rr}|\psi|^2~d x, \nonumber\\
 && I_3=\int\limits_{\rr}\left [ w^2+i(\phi^*\phi_x-\phi\phi^*_x
        +\psi^*\psi_x-\psi\psi^*_x)\right ]~d x, \nonumber\\
 && I_4=\int\limits_{\rr} \left [|\phi_x|^2+|\psi_x|^2+
          \beta(|\phi|^2+|\psi|^2)w \right ]~d x, \label{invariants}
 \eea
 where $I_1$ and $I_2$ are the mass functionals, $I_3$ is the momentum functional and $I_4$ is the energy
 functional, i.e. the Hamiltonian.  As the natural energy space for the LSI system
 (\ref{lwsw}) is $ H^1(\rr)\times H^1(\rr)\times L^2(\rr)$, the initial data  is taken
 as $(\phi_0,\psi_0,w_0) \in  H^1(\rr)\times H^1(\rr)\times L^2(\rr)$ in the present study.

 \setcounter{equation}{0}
 \section{Variational Characterization of Solitary Waves}

 \noindent
 In this section we briefly discuss a variational characterization of solutions
 for (\ref{cs1}), considered in \cite{borluk1}, where  the stability  analysis
 of solitary waves  is based on.

  Motivated by Nagy inequality \cite{nagy} given as
 \be
 \left( \frac{s}{2} H(\frac{s}{\beta}, \frac{p-1}{p}) \right )^{-\frac{\beta}{s}}\leq
         \frac{\| u_x \|_p ^\frac{\beta}{s}~ \| u \|_q ^{q+\beta\frac{q(p-1)}{ps}} }
           { \| u\|_{q+\beta}^{q+\beta}},  ~~~u\in H^1(\rr), \label{nagy}
 \ee
 where $q,\beta>0,~p\geq 1$, $\displaystyle s=1+q(p-1)/p$,
 $\displaystyle
 H(a,b)= [ (a+b)^{-(a+b)} \Gamma(1+a+b) ]/  [a^{-a} b^{-b} \Gamma(1+a)\Gamma(1+b) ]$
 and $\Gamma$ is the Gamma function, and by Gagliardo-Nirenberg inequality
 \bd
 \|u\|_r \leq C ~\|\nabla u\|_2^\vartheta ~ \|u\|_2^{1-\vartheta}, ~~~0<\vartheta \le 1, ~~~u\in H^1(\rr^n)
 \ed
 where $\displaystyle \vartheta=n(1/2-1/r)$; the  nonlinear functional
  $J(u,v)$  on $H^1(\rr) \times H^1(\rr)$
 \be
 J(u,v)= \frac{( \|u \|_2^2+ \|v \|_2^2)^{1-\theta/ 2}
          (\| u_x \|_2^2+ \| v_x \|_2^2)^{\theta /2}}
        { \|u^2+v^2\|_2^{1/2}  },~~~ \theta=\frac{1}{4},     \label{func}
 \ee
 was considered in \cite{borluk1}. The functional  $J(u,v)$ is well defined on
 $H^1(\rr) \times H^1(\rr)$ due to  embedding of $H^1(\rr)$ in $L^4(\rr)$.
 It should be pointed out that  the nonlinear  functional $J(u,v)$ is a generalization of
 the single variable functional $J(u)$ which was considered in the study of standing
 waves of the nonlinear Schr\"odinger equation \cite{wein3}.

 The first variation of the nonlinear functional  $J(u,v)$ is given as
 \bd
 \delta J=-B  \int\limits_\rr  \left\{
         [u_{xx}-\Omega  u+\gamma (u^2+ v^2) u]  \eta_1 
       +[v_{xx}-\Omega  v+\gamma (u^2+ v^2) v]  \eta_2  \right\}dx,  
 \ed
 where $\eta_i\in C_0^\infty(\rr)~(i=1,2)$, $~\displaystyle \Omega =\omega -c^2/4 $ and
 $~\displaystyle \gamma =\beta^2/c$,~
 $B= [ 3^3/(4^4 \Omega^3\gamma^4(\int\limits_\rr (u^2+ v^2)^2 dx)^6))]^{1/8}$,
 and the Pohozaev type identities,
 \be
   3 \int\limits_{\rr} ( u_x^2+v_x^2) dx = \Omega \int\limits_{\rr} (u^2+ v^2) d x
        =\frac{3\gamma}{4} \int\limits_{\rr}(u^2+v^2)^2   d x, \label{poh}
 \ee
 satisfied by $(u,v)$ are used. It is shown in \cite{borluk1} that
 the infimum of $J(u,v)$
 is achieved  at a  pair of positive functions $(R_1, R_2)$ when $c>0$ and $4\omega -c^2>0$
 using Lieb's compactness lemma.
 Thus the critical points of the functional $J(u,v)$ in $H^1(\rr)\times H^1(\rr)$ are
 the non-trivial weak solutions of (\ref{cs1}). As will be seen in section 4, the variational
 characterization plays a key role in the stability analysis of solitary waves (\ref{solution}).

 It should be noted that there are various  studies in the literature devoted to the problem
 of existence of solutions of the coupled system (\ref{cs1}) and its generalizations
 (\cite{maia, figu} and the references therein). In those studies, variational approaches
 based on minimization of  energy functionals subject to some constraints  are used.
 Though the approach followed in \cite{borluk1} is different from those of \cite{maia, figu},
 it is readily seen  that minimizing the energy functional is equivalent to minimizing the
 nonlinear functional $J(u,v)$. Indeed, the energy functional for solitary waves
  \bd
 I_4(u,v)=\int\limits_\rr \left ( u_x^2+v_x^2 +\frac{c^2}{4} (u^2+v^2)
    -\gamma  (u^2+v^2)^2\right ) dx,
 \ed
 after the scale transformation $(u_q(x), v_q(x))=\sqrt q  (u( q x), v( q x) )$ with
 $q > 0$, takes the form
 \bea
 &&\hspace{-0.5cm} I_4(u,v)\ge \inf\limits_{q >0} I_4(u_q, v_q) 
        =\inf\limits_{q >0}\int\limits_\rr [q^2( u_x^2+v_x^2) +\frac{c^2}{4} (u^2+v^2)
        -\gamma q(u^2+v^2)^2~] dx, \nonumber \\
        &&~~~~~~ \ge  \int\limits_\rr \left ( q^2( u_x^2+v_x^2)
        -\gamma q(u^2+v^2)^2\right ) dx, \label{ham0}
 \eea
 where the conserved mass integrals do not change: $\|u_q \|_2=\|u\|_2$ and
 $\|v_q\|_2=\|v\|_2$.  Using the scaled forms of the
 identities (\ref{poh})  in (\ref{ham0}), the energy functional
 takes the form
 \bd
   I_4(u,v)\ge \inf\limits_{q >0} I_4(u_q, v_q) \ge
       -\left (\frac{3 \gamma^2\Omega^7}{16}  \right )^{1/8}
       \lambda^{5/4}  \frac{1}{\inf J(u,v)},
 \ed
 for which $J(u_q, v_q)=J(u,v) $ and $\lambda =I_1+I_2$.
 Thus ground state solutions $(u_q, v_q)$, i.e. a minimizer of the Hamiltonian $I_4$,
 is  also  a minimizer of the functional $J(u,v)$.

 \setcounter{equation}{0}
 \section{Stability of Solitary Waves}

 \noindent
 In this section, we are concerned with  the stability of solitary wave solutions
 (\ref{solution}) of system  (\ref{lwsw}).
 For solitary waves, the appropriate notion of stability is orbital stability. All
 solitary waves of the same form but in different positions through space translation and
 phase rotation are assumed to be in the same orbit.  The LSI equations have
 translation  and phase  symmetries, i.e. if
 $~\left (\phi(x,t),\psi(x,t),w(x,t) \right )~$ solves the LSI equations, then
 $~  ( e^{i\theta_1}\phi(x+x_0,t), ~e^{i\theta_2}  \psi(x+x_0,t), ~w(x+x_0,t) )~$
 solves the same system for any $~x_0\in \rr$ and $\theta_1, \theta_2 \in [0,\,2\pi)$.
 We define the orbit ${\cal O}(f, g, h)$ of the triplet $(f,g,h)$ as follows:
 \bd
  {\cal O}(f, g,h)=\{ e^{i\theta_1} f(.+x_0), e^{i\theta_2} g(.+x_0), h(.+x_0);
   ~\theta_1,\theta_2 \in [0,\,2\pi), ~x_0\in \rr \}.
 \ed
 A solitary wave is said to be \emph{orbitally stable} if, for the initial data being
 near the solitary wave orbit, the solution at all later
 times remains near the solitary wave orbit.

  The main result of this section is the following theorem.

 \noindent
\textbf{Theorem 1.}
For $c>0$ and $4 \omega-c^2>0 $, solitary wave solution of
(\ref{lwsw})
  \bea
 \left . \begin{array}{ll}
 &   e^{i \omega t}\Phi(x-ct)=e^{i \omega t}R_1(x-ct) e^{i\frac{c (x-ct)}{2}},\\
 &   e^{i \omega t}\Psi(x-ct)=e^{i \omega t}R_2(x-ct) e^{i\frac{c (x-ct)}{2}},\\
 &\displaystyle    W(x-ct) = -\frac{\beta}{c}  [ R_1^2(x-ct)+ R_2^2(x-ct)],
 \end{array} \right \}  \label{solitary}
 \eea
 is orbitally stable,  i.e. for any $\epsilon\geq 0$ there exists a corresponding
 $\delta\geq 0$ such that the initial data
 $(\phi_0,\psi_0,w_0)\in H^1(\rr)\times H^1(\rr)\times L^2(\rr)$ with
 \bd
    \|\phi_0(.)-\Phi(.) \|_{H^1}\leq \delta,~~
    \|\psi_0(.)-\Psi(.) \|_{H^1}\leq \delta,~~
    \|w_0(.)-W(.)\|_2 \leq \delta,~~
 \ed
 imply
 \beas
 && \inf_{\substack{x_0 \in \rr \\ \theta_1 \in [0,2\pi)}}
        \|e^{i\theta_1}\phi(.+x_0,t)-\Phi(.)\|_{H^1} \leq \epsilon,\\
 && \inf_{\substack{x_0 \in \rr \\ \theta_2 \in [0,2\pi)}}
        \|e^{i\theta_2}\psi(.+x_0,t)-\Psi(.)\|_{H^1} \leq \epsilon,\\
 && \inf_{x_0 \in \rr}\|w(.+x_0,t)-W(.)\|_2 \leq \epsilon.
 \eeas

 In order to show that solitary waves (\ref{solitary}) are orbitally stable, i.e. to
 prove Theorem 1; we have to find an estimate on the distance in $H^1(\rr)\times H^1(\rr)$
 between the orbit  ${\cal O}_{(R_1,R_2)}$ of  solitary waves  and the solution
 $(\phi(x,t),\psi(x,t))$ of the LSI system.
 The deviation  of the solution $~(\phi(x,t), \psi(x,t))~$ corresponding to
 the initial data $(\phi_0,\psi_0)$  from the orbit ${\cal O}_{(R_1,R_2)}$ of
 solitary waves is measured by the metric
 \bd
 \rho_\Omega^2[(\phi, \psi),{\cal O}_{(R_1,R_2)}] =
    \inf_{\substack{x_0 \in \rr \\ \theta_1,\theta_2\in [0,2\pi)}}
    \left\{I_\Omega\right\},
 \ed
 where
 \bea
  I_\Omega(x_0,\theta_1,\theta_2) &=&
      N_\Omega (e^{i\theta_1} e^{-i\frac{c}{2}(.+x_0-ct)} \phi(.+x_0,t)-R_1) \nonumber \\
  &+& N_\Omega (e^{i\theta_2} e^{-i\frac{c}{2}(.+x_0-ct)} \psi(.+x_0,t)-R_2). \label{metric1}
 \eea
 The norm function $N_\Omega $  in (\ref{metric1}) is defined as
 $ N_\Omega(f)=\Omega\| f \|^2_2+ \| \nabla f \|^2_2$ and
 satisfies $  ~{\min (1,\Omega)} ~\| f \|^2_{H^1} \leq  N_\Omega(f) \leq
 {\max (1,\Omega)} ~\| f \|^2_{H^1}$. Perturbations of solitary waves, denoted by
 $~w_1(x,t), ~w_2(x,t)$ and $\eta(x,t)$,  are defined in the form
 \bea
 && w_1(x,t) = e^{i\theta_1}e^{-i\frac{c}{2}(x+x_0-ct)} \phi(x+x_0,t)-R_1(x),  \label{sol1} \\
 && w_2(x,t) = e^{i\theta_2}e^{-i\frac{c}{2}(x+x_0-ct)} \psi(x+x_0,t)-R_2(x),  \label{sol2}\\
 &&  \eta(x,t)=u(x+x_0,t)+\frac{\beta}{c} \left [ R_1^2(x)+ R_2^2(x) \right ], \label{sol3}
 \eea
 where $~w_k(x,t)=  p_k(x,t)+iq_k(x,t)~(k=1,2)$ are complex-valued functions, and
 $~\eta(x,t)~$ is a real-valued function. Here  $\theta_1, \theta_2$ and $x_0$ will be chosen
 later where the infimum of  $I_\Omega$ is attained. Eq. (\ref{metric1}), and  (\ref{sol1})-(\ref{sol3})
 show that we have  to find estimates on the  $H^1$ norms of $w_1(x,t)$ and $w_2(x,t)$, and the
 $L^2$ norm of  $\eta(x,t)$.

 The following lemma  is a generalization of the one which was proved in the context
 of the orbital stability of solitary waves, by Bona \cite{bona} for  the Korteweg-de Vries
 equation and by Angulo and Montenegro \cite{pava1}  for the long wave-short wave
 interaction equations with an integral term. The following lemma states that there are
 $\theta_i=\theta_i(t)~(i=1,2)$ and $x_0=x_0(t)$ such that infimum of
 $I_\Omega(x_0,\theta_1,\theta_2)$ exists where the local well-posedness of the Cauchy
 problem for (\ref{lwsw}) is used.

 \noindent
 \textbf{ Lemma 2.} Let $(\phi,\psi,u)$ be a solution of (\ref{lwsw}) corresponding to
 the initial data
 $(\phi_0,\psi_0,u_0) \in H^1(\rr)\times H^1(\rr)\times L^2(\rr)$ with the properties
 $  \|\phi_0\|_2= \|R_1 \|_2~$ and $~\|\psi_0\|_2= \|R_2 \|_2 $.
 Suppose that $~ I_\Omega(x_0,\theta_1,\theta_2)<\Omega(\|R_1\|_2^2+\|R_2\|_2^2)~$
 for some $t_0 \in [0,T]$ and some
 $(x_0,\theta_1,\theta_2)\in\rr\times [0,2\pi)\times[0,2\pi)$. Then
 $~\mbox{inf}\{I_\Omega |x_0\in \rr, \theta_1, \theta_2\in [0,2\pi)\} ~$
 is assumed at least once.

\textbf{Proof.} It is clear that $I_\Omega $ is a continuous function of
 $(x_0,\theta_1, \theta_2)$ on $\rr \times [0,2\pi) \times [0,2\pi) $.
 Moreover, for any $(\theta_1,\theta_2)\in [0,2\pi)\times[0,2\pi)$, we have
 \bea
 \lim_{x_0\rightarrow\mp \infty}I_\Omega(x_0,\theta_1,\theta_2)&=&
    \|[e^{-i\frac{c}{2}(\cdot-ct)}\phi(\cdot,t)]'\|_2^2+
   +\|[e^{-i\frac{c}{2}(\cdot-ct)}\psi(\cdot,t)]'\|_2^2 \nonumber \\
    &&+\|R_1'(\cdot)\|_2^2+\|R_2'(\cdot)\|_2^2
      +2\Omega\|R_1(\cdot)\|_2^2+2\Omega\|R_2(\cdot)\|_2^2, \nonumber \\
  &&=(\frac{7 \Omega}{3}+\frac{c^2}{4})(\|R_1\|_2^2+ \|R_2\|_2^2), \label{iomega}
 \eea
 where (\ref{poh}) is used.
 The hypothesis $~ I_\Omega(x_0,\theta_1,\theta_2)<\Omega(\|R_1\|_2^2+\|R_2\|_2^2)~$,
 the continuity of $I_\Omega$  and (\ref{iomega}) imply the result.$\hfill\square$

 We now  show that the infimum of  $I_\Omega$ is attained at a finite value of $x_0$ for
 some $t_0 \in [0,T]$. For this aim, it will suffice to show that
 $~ I_\Omega(x_0,\theta_1,\theta_2)<\Omega(\|R_1\|_2^2+\|R_2\|_2^2)~$ holds  in some interval.
 Using the inequality $\|a+b\|_2^2\leq 2\|a\|_2^2+2\|b\|_2^2$,
 one can obtain
 \beas
 I_\Omega(ct, -\omega t, -\omega t)&\leq &
        2\|\phi'(\cdot)-\phi_s'(\cdot)\|_2^2+(\frac{c^2}{2}+\Omega)\|\phi(\cdot)
        -\phi_s(\cdot)\|_2^2 \nonumber\\
 &&
        +2\|\psi'(\cdot)-\psi_s'(\cdot)\|_2^2+(\frac{c^2}{2}+\Omega)\|\psi(\cdot)
        -\psi_s(\cdot)\|_2^2, 
 \eeas
 where  prime denotes differentiation with respect to spatial variable $x$.
 Solitary wave solutions $~(\phi_s,\psi_s)~$ given in (\ref{solution})
 are globally defined. Thus, it follows from the continuous dependence theory that,
 for a $T>0$, there exists a $\delta>0$ such that if
 $$
 \|\phi_0(.) - e^{i\frac{c}{2}\cdot} R_1(.) \|_{H^1}< \delta,~~~ \mbox{and}~~~
 \|\psi_0(.) - e^{i\frac{c}{2}\cdot} R_2(.)\|_{H^1}< \delta,
 $$
 then the solution $(\phi(x,t),\psi(x,t))$ corresponding to the initial data
 $(\phi_0(x),\psi_0(x))$ exists at least for $0\leq t\leq T$. This solution also satisfies
 $$
 \| \phi(\cdot,t) -\phi_s(\cdot,t)\|_{H^1}< \epsilon,~~~\mbox{and}~~~
 \| \psi(\cdot,t) -\psi_s(\cdot,t)\|_{H^1}< \epsilon.
 $$
 Using this result, we get $ I_\Omega(ct, -\omega t, -\omega t)\leq 4\epsilon^2(1+\omega).$
 Choosing $ \epsilon^2<\Omega(\|R_1\|_2^2+\|R_2\|_2^2)/[4(1+\omega)],$
 shows that the hypothesis of  Lemma 2 is satisfied at least for
 $(\tilde{x}_0,\tilde{\theta}_1,\tilde{\theta}_2) =(ct, -\omega t, -\omega t)$,
 from which we get an upper bound for  $I_\Omega$.

 As a result of Lemma 2, the following compatibility conditions are
 obtained  for the real-valued increment functions $p_i(x,t)$ and $q_i(x,t)~~(i=1,2)$
 \bea
 && \int\limits_\rr  \left ( R_1^2+ R_2^2 \right )R_1 q_1dx=0, \label{eq31}\\
 && \int\limits_\rr  \left ( R_1^2+ R_2^2 \right )R_2 q_2dx=0, \label{eq32} \\
 && \int\limits_\rr (R_1^2+ R_2^2)\left ( R_1 \frac{\partial p_1}{\partial x}
    +R_2\frac{\partial p_2}{\partial x} \right )dx=0.\label{eq33}
 \eea
 The relations  (\ref{eq31}),  (\ref{eq32})  and (\ref{eq33}) are obtained by differentiating
 $I_\Omega$  defined in (\ref{metric1}) with respect to  $\theta_1,\theta_2$ and $x_0$,
 using system (\ref{cs1}) and  then evaluating  the resulting equations at values
 $(x_0, \theta_1, \theta_2)$ which  minimize  $I_\Omega$. Note that
 \bd
 I_\Omega=  \|e^{i\theta_1} A'-R_{1}' \|_2^2+ \|e^{i\theta_1} B'-R_{2}'\|_2^2
 +\Omega  \|e^{i\theta_1} A-R_1\|_2^2  +\Omega  \|e^{i\theta_1} B-R_2\|_2^2,
 \ed
 where  $e^{i \theta_1} e^{-i\frac{c}{2}(x+x_0-ct)}\phi(x+x_0,t)= e^{i \theta_1}A(x+x_0,t)=R_1(x)+w_1(x,t)~$ and
        $~e^{i \theta_2} e^{-i\frac{c}{2}(x+x_0-ct)}\psi(x+x_0,t)= e^{i \theta_2}B(x+x_0,t)=R_2(x)+w_2(x,t)$.

 We now introduce a continuous nonlinear  functional $L$, called the Lyapunov functional,
 over  $H^1(\rr)\times H^1(\rr)\times L^2(\rr)$ in the form
 \be
 L \left (\phi,\psi,u \right )=\omega ( I_1 + I_2) +\frac {c}{2} I_3+I_4, \label{lyapunov}
 \ee
 where $I_k~(k=1,2,3,4)$, given in (\ref{invariants}), are the conserved quantities of
 system (\ref{lwsw}). Thus,  the Lyapunov functional is invariant with
 time: $\Delta L(0)= \Delta L(t)$.  Our stability result will rely on the inequalities
 \beas
 && \Delta L(0)\leq 2 g(\varepsilon)\\
 && \Delta L(t)\ge g(\|w_1\|_{H^1})+ g(\|w_2\|_{H^1}),
 \eeas
 where $g(x)=a_1 x^2-a_2 x^3-a_3 x^4$ for some positive constants $a_i~(i=1,2,3)$, and
 $\|w_i\|_{H^1}~(i=1,2)$ is the distance between the solitary wave $(\Phi, \Psi)$ and the
 solution $(\phi, \psi)$  of (\ref{lwsw}).
 To find the bounds, we calculate $\Delta L(t)$
  \beas
  \Delta L(t)&=& L(\phi(x,t),\psi(x,t),u(x,t))-L(\Phi(x),\Psi(x), U(x)),\\
    &=& L\left (\Phi(x)+e^{i\frac{c x}{2}}w_1(x,t),\Psi(x)+e^{i\frac{c x}{2}}w_2(x,t), U(x)+\eta(x,t))\right.\\
     && \!\!\!\!\!\!\!\!\!\!\left. -L(\Phi(x),\Psi(x), U(x)\right ).
 \eeas
 Expanding the functional $L$ near $(\Phi,\Psi)$ yields
 \be
 \Delta L(t)= \delta L+ \delta^2 L+ \delta^3 L, \label{delly}
 \ee
 where $~ \delta L,~ \delta^2 L$ and
 $ \delta^3 L$ are the first, second and third variations of $L$,
 respectively; and  all variations higher than third order are zero.
 The explicit forms of variations are given as
 \bea
  \delta L &=& \int\limits_\rr 2 \left \{
         \left [ R_{1,xx}-\Omega R_1 + \gamma (R_1^2+R_2^2) R_1\right ] p_1 \right.\nonumber\\
       &&  \!\!\!\!\!\!\!\!\!\! +\left. \left [ R_{2,xx}-\Omega R_2 + \gamma (R_1^2+R_2^2) R_2\right ] p_2
        \right \} dx, \label{var1} \\
  \delta^2 L &=& \int\limits_\rr
        \left [\frac{c}{2} \eta^2 + p_{1,x}^2+q_{1,x}^2 + p_{2,x}^2+q_{2,x}^2
        + \Omega (p_1^2+q_1^2+ p_2^2 + q_2^2)\right. \nonumber\\
      &&\!\!\!\!\!\!\!\!\!\! \left.  +2\beta (R_1 p_1 + R_2 p_2) \eta
   -\gamma (R_1^2 + R_2^2)(p_1^2+q_1^2+p_2^2+q_2^2) \right ]dx,
        \label{var2}\\
  \delta^3 L &=& \int\limits_\rr
         \beta (p_1^2+q_1^2 +p_2^2+q_2^2) \eta~ dx, \label{var3}
 \eea
 where the relations $\displaystyle
 \Phi(x)=R_1(x)e^\frac{icx}{2},~ \Psi(x)=R_2(x)e^\frac{icx}{2},~ W(x)=-\beta^2(R_1^2(x)+ R_2^2(x))/c$, and
 $w_k(x)=p_k(x)+ i q_k(x)~(k=1,2)$ are used.  Because  $(R_1,R_2) $ is  a solution (\ref{cs1}),
 the first variation (\ref{var1}) vanishes. Thus $(R_1,R_2) $  is also a critical point of the
 Lyapunov functional $L$. From eqs.  (\ref{var2}) and (\ref{var3}), we have
 \bea
 \Delta L(t)&=& \langle L_0 q_1,q_1 \rangle+\langle L_0 q_2,q_2 \rangle
         +\langle L_1 p_1,p_1 \rangle + \langle L_2 p_2,p_2 \rangle
        + 2 \langle L_3 p_1,p_2 \rangle   \nonumber\\
 &&\!\!\!\!\!\!-\gamma \int\limits_\rr \left [ \frac{1}{2}(p_1^2+q_1^2+ p_2^2+q_2^2)^2
  + 2 (p_1^2+q_1^2+p_2^2+q_2^2)(p_1 R_1+p_2 R_2) \right ]dx \nonumber\\
 &&\!\!\!\!\!\!+\frac{c}{2}\int\limits_\rr  \left [ \eta+\frac{2 \beta}{c} (p_1 R_1+ p_2R_2)
        +\frac{\beta}{c}(p_1^2+q_1^2+p_2^2+q_2^2) \right ]^2dx,~ \label{lvar}
 \eea
 where  the operators $L_i~(i=0, 1,2,3)$ are defined as
 \beas
  && L_0=-\frac{\partial^2}{\partial x^2}+\Omega -\gamma (R_1^2 +
  R_2^2),~~~~
   L_1=-\frac{\partial^2}{\partial x^2}+\Omega -\gamma (3 R_1^2 +
   R_2^2),\\
  && L_2=-\frac{\partial^2}{\partial x^2}+\Omega -\gamma ( R_1^2 +
  3R_2^2),~~~~
   L_3=-2\gamma R_1 R_2.
 \eeas
 We use the following lemmas to find a lower bound for $\Delta L(t)$.

 \noindent
 \textbf{Lemma 3.} There exist positive constants $C_i~(i=1,2)$ such that
 \be
 \langle L_0 q_i, q_i \rangle ~\geq  C_i \|q_i \|_{H^1}^2 ~~~(i=1,2).      \label{qiqi}
 \ee

 \textbf{Proof.} It should be noted that
 $L_0 R_i=0$ and  $~R_i>0 ~ (i=1,2)$. Therefore $L_0$ is
 a non-negative operator, i.e.
 $~\displaystyle \mu_i =\inf ( \langle L_0 q_i, q_i \rangle / \langle q_i,q_i\rangle )\geq 0~ (i=1,2)$.
 If the infimum  of the functional $~ \mu_i$ subject to the
 constraints (\ref{eq31}) and (\ref{eq32})
 is zero then it is attained at $ q_i(x)=R_i(x)$.
 This contradicts to the above constraints, thus $\mu_i>0~(i=1,2)$, i.e.
 \bd
 \langle L_0 q_i, q_i \rangle =\frac{1}{k_i+1} \| q_i \|-\gamma \int\limits_\rr (R_1^2+R_2^2)q_1^2 dx
    +\frac{k_i}{k_i+1}\| q_i \| ~\geq  \bar C_i \| q_i \|_2^2 ~~~(i=1,2)
 \ed
 where $\|q_i\|= \|\nabla q_i \|_2^2+\Omega \|q_i\|_2^2$, $k_i $ and $\bar C_i$ are some positive constants.
 If $ ~\displaystyle k_i< {\bar C_i} /(2\gamma E^2)~$ where $~E=\max(\|R_1\|_\infty,\|R_2\|_\infty)$,
 then we have  $~  \| q_i \| /(k_i+1)-\gamma \int\limits_\rr (R_1^2+R_2^2)q_1^2 dx >0$,
 and consequently
 \bd
 \langle L_0 q_i, q_i \rangle ~\geq  C_i  \|q_i \|_{H^1}^2 ~~~(i=1,2),
 \ed
 where $~\displaystyle C_i= k_i \min(1,\Omega)/ (k_i+1)$.
$\hfill\square$

 \noindent
 To find a lower bound for the expression $\langle L_1 p_1, p_1 \rangle
 +\langle L_2 p_2, p_2 \rangle+2\langle L_3 p_1, p_2 \rangle$ in   (\ref{lvar})
 is more difficult than that of $\langle L_0 q_i, q_i \rangle$. We will use the facts
 that $(R_1, R_2)$ is the minimizer of the functional $J(u,v)$ and that the
 expression $\langle L_1 p_1, p_1 \rangle +\langle L_2 p_2, p_2 \rangle
 + 2\langle L_3 p_1, p_2 \rangle$ is associated with the second variation of $J(u,v)$.
 First we prove the following lemma which is a generalization of the one given
 in \cite{wein2}.

 \noindent
\textbf{Lemma 4.}  $ \displaystyle \inf_{\substack{{\langle f,R_1\rangle=0}\\
        {\langle g,R_2\rangle=0}}} \left ( \langle L_1 f,f \rangle
 +\langle L_2 g,g \rangle+2 \langle L_3 f,g \rangle\right )=0$.

 \textbf{Proof.} Recall that $(R_1, R_2)$ is a minimizer of the
 nonlinear functional $J(u,v)$. Thus $~\delta^2 J\geq 0$
 near $(R_1,R_2)$. The second variation of the functional $J$ is of the form
 \bea
  \frac{d^2}{d \epsilon^2}
     J(R_1+\epsilon \eta_1,\,R_2+\epsilon\eta_2) \mid_{\epsilon=0}
     &=& a^2 \left ( \langle L_1  \eta_1,\eta_1\rangle
     +      \langle L_2 \eta_2 ,\eta_2\rangle
     + 2    \langle L_3 \eta_1,\eta_2  \rangle \right ) \nonumber \\
  && \hskip-4cm
  +  a^2   \left [ \frac{\Omega^2}{3 d} (\langle R_1,\,\eta_1\rangle + \langle R_2,\,\eta_2\rangle)^2 \right .
       \nonumber \\
  && \hskip-4cm
  \left. + \frac{2\Omega}{d} (\langle R_1,\,\eta_1\rangle + \langle R_2,\,\eta_2\rangle)
     (\langle R_{1,x},\,\eta_{1,x} \rangle
     + \langle R_{2,x},\,\eta_{2,x}\rangle )     \right ] \nonumber \\
  && \hskip-4cm
  - \frac{a^2}{d}  \left( \langle R_{1,x},\,\eta_{1,x} \rangle
        + \langle R_{2,x},\,\eta_{2,x}\rangle  \right )^2     \geq 0, \label{del}
 \eea
 where $~\displaystyle a^2= [ 27 \gamma^2/ (\Omega^3 d^6)]^{1/8}/(4\sqrt 2) $,
 and $d=  \int\limits_\rr (u_x^2+v_x^2) dx$.
 It should be noted that eq. (\ref{cs1}) and Pohozaev type identities given by
 (\ref{poh}) are used in obtaining (\ref{del}).

 If the increment functions  are chosen as $\eta_1=f$ and $\eta_2=g$ with the properties
 $~\langle f,R_1\rangle =0$ and $\langle g,R_2\rangle=0$, then it follows from (\ref{del}) that
 \be
    \langle L_1 f,f\rangle +\langle L_2 g,g \rangle
        + 2 \langle L_3 f,g \rangle ~\geq 0.\label{exp}
 \ee
 Moreover,  the functions $R_{1,x}$ and $R_{2,x}$ satisfy
  \be
 \left . \begin{array}{ll}
 & L_1  R_{1,x} +L_3  R_{2,x}=
    \left ( -R_{1,xx}+\Omega R_1- \gamma (R_1^2 +R_2^2) R_1 \right) _x=0, \\
 & L_2  R_{2,x} +L_3  R_{1,x}=
    \left ( -R_{2,xx}+\Omega R_2- \gamma (R_1^2 +R_2^2) R_2 \right) _x=0.
 \end{array} \right \} \label{l1l3}
 \ee
 As a result of  (\ref{l1l3}) we find
 \beas
 &&  \langle L_1  R_{1,x},R_{1,x} \rangle+ \langle L_2 R_{2,x},R_{2,x}\rangle
    +\langle L_3  R_{1,x},R_{2,x} \rangle+ \langle L_3 R_{2,x},R_{1,x}\rangle \nonumber \\
 && ~~~~~~~~~ =\langle L_1 R_{1,x}+L_3 R_{2,x},R_{1,x}\rangle
              +\langle L_2 R_{2,x}+L_3 R_{1,x},R_{2,x}\rangle=0,
 \eeas
 which shows that the infimum of (\ref{exp}) is assumed at $(R_{1,x}, R_{2,x})$.
 Because $f=R_{1,x}$  and $g=R_{2,x}$ satisfy the hypothesis of the lemma, we get
 $\langle f,R_1 \rangle= \langle R_{1,x}, R_1\rangle=0$
 and $\langle g,R_2\rangle=\langle R_{2,x}, R_2\rangle=0$. This completes the
 proof.
$\hfill\square$

 In order to find a  lower bound for
 $\langle L_1  p_1,p_1\rangle + \langle L_2  p_2, p_2\rangle+ 2\langle L_3  p_1, p_2\rangle$,
 we  require that  the perturbed solution has the same $L^2$-norm as the
 solitary wave, as given in the hypotheses of Lemma 2,
 \be
    \| \phi\|_2= \| R_1\|_2,~~~\| \psi\|_2= \| R_2\|_2. \label{res}
 \ee
 Conditions (\ref{res}) give rise to the following constraints
 \be
 \langle R_i,p_i \rangle=-\frac{1}{2}[\langle p_i,p_i\rangle+ \langle q_i,q_i\rangle]
    =-\frac{1}{2} \| w_i \|_2^2 ~<0~~~(i=1,2), \label{cond}
 \ee
 where definitions (\ref{sol1}) are used. The restrictions  (\ref{res}) will be relaxed later
 and the stability of solitary waves will be proved with respect to general perturbations.
 To this end, we assume that the real parts of the increment functions,
 $p_i(x,t)~(i=1,2)$, will be of the form $~ p_i=p_{i|\!|}+p_{i\bot}~$ where
 \bd
 p_{i|\!|}=     \frac{\langle p_i,\,R_i \rangle}{\|R_i\|^2_2} R_i,~~~
 p_{i\bot}= p_i-\frac{\langle p_i,\,R_i \rangle}{\|R_i\|^2_2} R_i.
 \ed
 This gives rise to $\langle p_{i\bot}, R_i \rangle=0~~(i=1,2)$. Using the decomposition of
 $p_i(x,t)~(i=1,2)$, we have
 \bea
&& \langle L_1 p_1,p_1\rangle+\langle L_2 p_2,p_2\rangle+2\langle L_3 p_1,p_2\rangle \nonumber\\
       &&~~~~~~~~~~~~~~~~~~= \langle L_1 p_{1\bot},p_{1\bot}\rangle
          +\langle L_2 p_{2\bot},p_{2\bot}\rangle
         +2\langle L_3 p_{1\bot},p_{2\bot}\rangle\nonumber\\
       && ~~~~~~~~~~~~~~~~~~+\langle L_1 p_{1|\!|},p_{1|\!|} \rangle
          +\langle L_2 p_{2|\!|},p_{2|\!|}\rangle
         +2\langle L_3 p_{1|\!|},p_{2|\!|} \rangle
         +2\langle L_1 p_{1\bot},p_{1|\!|} \rangle\nonumber\\
       &&~~~~~~~~~~~~~~~~~~+2\langle L_2 p_{2\bot},p_{2|\!|}\rangle
         +2\langle L_3 p_{2|\!|},p_{1\bot} \rangle
         +2\langle L_3 p_{1|\!|},p_{2\bot}\rangle. \label{Li+}
 \eea
 To find a suitable lower bound for
 $\langle L_1 p_{1\bot},p_{1\bot}\rangle +\langle L_2 p_{2\bot},p_{2\bot}\rangle
       +2\langle L_3 p_{1\bot},p_{2\bot}\rangle$ using Lemma 4, we further  assume that
 $~\langle p_1,\,R_1 \rangle/ \|R_1\|^2_2= \langle p_2,\,R_2 \rangle /\|R_2\|^2_2$.

 \noindent
\textbf{Lemma 5.} There exist  positive constants $C_3$ and $C_4$  such that
 \bea
 && \langle L_1 p_{1\bot},p_{1\bot}\rangle+\langle L_2 p_{2\bot},p_{2\bot}\rangle
    +2 \langle L_3 p_{1\bot},p_{2\bot} \rangle \nonumber \\
 &&~~~~~~~~~~    \geq C_3 ( \|p_1 \|_2^2  +\|p_2\|_2^2)
    - C_4( \|w_1\|_{H^1}^4 - \|w_2\|_{H^1}^4). \label{ineq1}
 \eea

\noindent
 \textbf{Proof.}
 If $f=p_{1\bot}$ and $g=p_{2\bot}$ then the hypotheses of Lemma 4 are satisfied by
 $p_{1\bot}$ and $p_{2\bot}$. That is,
 \be
 \langle L_1 p_{1\bot},p_{1\bot}\rangle+\langle L_2 p_{2\bot},p_{2\bot}\rangle
        +2 \langle L_3 p_{1\bot},p_{2\bot} \rangle \geq 0.    \label{pdpd}
 \ee
 The infimum of  (\ref{pdpd}) is zero. 
 This infimum is attained at  $(p_{1\bot},p_{2\bot})=(R_{1,x},
 R_{2,x})$.  In such a case,  for the increment functions
 $p_i=\alpha R_i+R_{i,x}~(i=1,2)$ where $\alpha=\langle p_i,R_i \rangle /\|R_i\|_2^2~(i=1,2)$,
 the constraint (\ref{eq33}) reduces to
 \beas
 &&  \frac{\alpha}{4} \int\limits_\rr [(R_1^2+ R_2^2)^2]_x dx
        + \int\limits_\rr (R_1^2+ R_2^2)(R_1 R_{1xx}+ R_2 R_{2xx}) dx=0,\\
 &&  \int\limits_\rr \left \{ [ (R_1^2+ R_2^2)_x]^2
    + 2 (R_1^2+ R_2^2) [(R_{1,x})^2+ (R_{2,x})^2]  \right \} dx=0,
 \eeas
 where integration by parts is used.
 This result leads to $R_i=0~(i=1,2)$ which  contradicts  positivity
 of ground state solutions  $(R_1,R_2)$.  Thus there exists a positive  constant $\bar C_3$ such that
 \be
 \langle L_1 p_{1\bot},p_{1\bot}\rangle+\langle L_2 p_{2\bot},p_{2\bot}\rangle
        +2 \langle L_3 p_{1\bot},p_{2\bot} \rangle \geq \bar C_3.\label{Ldik1}
 \ee
 Moreover, using  $ \displaystyle \langle p_{i\bot},p_{i\bot} \rangle = \langle p_i,p_i\rangle
 -[\langle  p_i,p_i \rangle + \langle  q_i,q_i \rangle ]^2/(4 \|R_i\|_2^2)$, the
 inequality (\ref{Ldik1}) can be arranged  to yield (\ref{ineq1})
 \beas
 && \langle L_1 p_{1\bot},p_{1\bot}\rangle+\langle L_2 p_{2\bot},p_{2\bot}\rangle
    +2 \langle L_3 p_{1\bot},p_{2\bot} \rangle
 \geq   C_3  (\langle p_{1\bot},p_{1\bot} \rangle +
     \langle p_{2\bot},p_{2\bot} \rangle  ), \nonumber \\
 && \hskip4cm= C_3 \left (\| p_1 \|_2^2 + \| p_2 \|_2^2
     -\frac{\|w_1\|_2^4}{4 \|R_1\|_2^2 } -\frac{\|w_2\|_2^4}{4 \|R_2\|_2^2 } \right ),\nonumber \\
 && \hskip4cm\geq  C_3  ( \|p_1 \|_2^2  +\|p_2\|_2^2) -C_4(\|w_1\|_{H^1}^4+ \|w_2\|_{H^1}^4),
 \eeas
 where continuous embedding of $H^1(\rr)$ in $L^4(\rr)$ is used, and
 $C_3$  and $C_4$ are some positive constants. This completes the proof of Lemma 5.
$\hfill\square$

 \noindent
 \textbf{Lemma 6.}  There exist  positive constants $C_5$  and  $C_6$  such that
 \bea
   \langle L_1 p_{1|\!|},p_{1|\!|}\rangle+ \langle L_2 p_{2|\!|},p_{2|\!|}\rangle
   +2 \langle L_3 p_{1|\!|},p_{2|\!|}\rangle
        \geq  -C_5 \|w_1\|_{H^1}^4 - C_6 \|w_2\|_{H^1}^4.  \label{ineq21}
 \eea

\noindent
 \textbf{Proof.} Recall that
 $\langle L_i R_i,R_i \rangle = -2 \gamma \langle R_i^2,R_i^2 \rangle ~(i=1,2)$.
 Firstly, using $  p_{i|\!|}=\alpha R_i  ~(i=1,2)$ where $~\alpha=- \|w_i\|_2^2/(2\|R_i\|_2^2)$,
 we obtain
 \be
  \langle L_i p_{i|\!|},p_{i|\!|} \rangle = \alpha^2 \langle L_i R_i,R_i  \rangle=
    -\frac{\gamma}{2} \frac{\|R_i^2\|_2^2}{\|R_i\|_2^4}~\|w_i\|_2^4
    \geq - \bar C_{4+i} \|w_i\|_{H^1}^4  ~~(i=1,2), \label{ineq2}
 \ee
 where $\bar C_5$ and $\bar C_6$ are positive constants.  Secondly,  using Sobolev embedding  and Young's
 inequality $~ ab\leq a^p/p+ b^q/ q~$ with $~p=q=2~$, we obtain
 \bea
  \langle L_3 p_{1|\!|},p_{2|\!|} \rangle
  =-\frac{\gamma}{4}~\frac{ \|R_1 R_2 \|_2^2}{\|R_1\|_2^2 ~\|R_2\|_2^2}~\|w_1\|_{2}^2 ~\|w_2\|_{2}^2 ~
  \geq -\frac{\bar C_7}{2} \left( \|w_1\|_{H^1}^4+\|w_2\|_{H^1}^4\right ), \label{ineq3}
 \eea
 where $\bar C_7$ is a positive constant. (\ref{ineq21}) follows from (\ref{ineq2}) and (\ref{ineq3}).
 $\hfill\square$

 \noindent
 \textbf{Lemma 7.} $  \langle L_3  p_{1|\!|},p_{2\bot} \rangle = 0$ and
 $\langle L_3  p_{2|\!|},p_{1\bot} \rangle=0$.

 \noindent
 \textbf{Proof.} Using the definition of the operator $L_3$, we have
 $ \langle L_3 p_{1|\!|},p_{2\bot} \rangle = -2\gamma \alpha \langle R_2^2 p_{1\bot}, R_1\rangle$.
 Then
 \be
  |\langle L_3 p_{1|\!|},p_{2\bot} \rangle| \leq |2\gamma \alpha| E^2 |\langle  p_{1\bot},R_1 \rangle|=0
    \label{ineq6}
 \ee
 and, similarly
 \be
 |\langle L_3 p_{2|\!|},p_{1\bot} \rangle|\leq
    |2\gamma \alpha| E^2 |\langle  p_{2\bot},R_2 \rangle|=0.    \label{ineq7}
 \ee
 This completes the proof.
 $\hfill\square$

 \noindent
 \textbf{Lemma 8.} There exist positive constants  $E_i$ and $F_i~(i=1,2)$ such that
  \bea
   2\langle L_i p_{i\bot} ,p_{i|\!|} \rangle
        \geq  -E_i \|w_i\|_{H^1}^3 -F_i \|w_i\|_{H^1}^4~~ (i=1,2).  \label{ineq31}
 \eea

 \noindent
 \textbf{Proof.}  For the terms $\langle L_i p_{i\bot},p_{i|\!|} \rangle$, we find
 \be
 \langle L_ip_{i\bot},p_{i|\!|} \rangle =\alpha \left (\langle R_{i,x}, p_{i\bot,x}\rangle
    -3 \gamma \langle R_i^3, p_{i\bot}\rangle - \gamma \langle R_j^2 R_i, p_{i\bot}\rangle \right ),
    ~~~(i,j=1,2~~i\neq j), \label{lem8}
 \ee
 where $\alpha=-\|w_i\|_2^2/(2 \|R_i\|_2^2)$,
 $~|\langle R_i^3, p_{i\bot}\rangle|\leq E^2 |\langle R_i, p_{i\bot}\rangle|=0$,  and
 $~|\langle R_j^2 R_i, p_{i\bot}\rangle|\leq E^2 |\langle R_i, p_{i\bot}\rangle|=0$.
 Using $p_{i\bot}=p_i-\alpha R_i$ and the Cauchy-Schwartz inequality in (\ref{lem8}), we have
 \beas
 \langle L_i p_{i\bot},p_{i|\!|} \rangle
 &\geq &  -\frac{ \|w_i\|_2}{2 \|R_i\|_2^2} \langle R_{i,x}, p_{i,x}\rangle
          -\frac{ \|R_{i,x}\|_2^2}{4 \|R_i\|_2^4}~\|w_i\|_2^4,  \\
 &\geq &  -\frac{ \|R_{i,x}\|_2}{2 \|R_i\|_2^2}~ \|w_i\|_2^2~ \|w_{i,x}\|_2
          -\frac{ \|R_{i,x}\|_2^2}{4 \|R_i\|_2^4 } ~\|w_i\|_2^4~~~(i=1,2).
 \eeas
 By  continuous embedding of $H^1(\rr)$ in $L^2(\rr)$ the result follows
 \bd
 \langle L_i p_{i\bot},p_{i|\!|} \rangle \geq
    -\frac{E_i}{2}\|w_i\|_{H^1}^3 -\frac{F_i}{2} \|w_i\|_{H^1}^4~~~~(i=1,2), 
 \ed
 where $E_i$ and $F_i$ are some positive constants .
$\hfill\square$

 \noindent
 \textbf{Lemma 9.} There exist positive constants  $A_i~~(i=1,2,3)$ such that
 \bea
 \langle L_1 p_1,p_1 \rangle
    +\langle L_2 p_2,p_2 \rangle+2\langle L_3 p_1,p_2 \rangle
 &\geq&
    A_1 \left( \| p_1 \|_{H^1}^2 +\|p_2\|_{H^1}^2 \right )\nonumber\\
 &&\hskip-5cm - A_2 \left( \| w_1 \|_{H^1}^3 +\|w_2\|_{H^1}^3 \right )
 -A_3 \left( \| w_1 \|_{H^1}^4 +\|w_2\|_{H^1}^4 \right). \label{L+}
 \eea

 \noindent
 \textbf{Proof.} By direct computation, one can see that
 \bea
 &&\langle L_1 p_1,p_1 \rangle+  \langle L_2 p_2,p_2\rangle
      +  2\langle L_3  p_1,p_2\rangle\nonumber \\
 &&~~~~~~~~~~ = -\gamma  \int\limits_\rr [(R_1^2+R_2^2) (p_1^2+ p_2^2) + 2 (R_1 P_1+ R_2 p_2)^2]  dx\nonumber \\
  &&~~~~~~~~~~   + ~\|p_1 \|+ \|p_2 \|, \label{eq8}
 \eea
 where $~  \|p_i \|= \Omega \|p_i \|_2^2+ \|\nabla p_i \|_2^2~~(i=1,2)$.
 On the other hand, combining the inequalities (\ref{ineq1}), (\ref{ineq21}),  (\ref{ineq6}),
  (\ref{ineq7}) and  (\ref{ineq31}),  we obtain
 \bea
 \langle L_1 p_1,p_1 \rangle+  \langle L_2 p_2,p_2\rangle
      +  2\langle L_3  p_1,p_2\rangle
  &\geq & C_3(\|p_1 \|_2^2+ \|p_2 \|_2^2)\nonumber \\
  && -E_1\| w_1 \|_{H^1}^3 - E_2\| w_2 \|_{H^1}^3\nonumber \\
  && -   C_8\| w_1 \|_{H^1}^4  -  C_9\| w_2 \|_{H^1}^4, \label{eq9}
 \eea
 where $C_8=C_4+C_5+F_1$ and $C_9=C_4+C_6+F_2$ are positive constants.

 Using (\ref{eq8}) and (\ref{eq9}), for a sufficiently small positive number $m$, we find
 \bea
 && I=\frac{1}{m+1}(\|p_1 \|+ \|p_2 \|)
    - \gamma  \int\limits_\rr [(R_1^2+R_2^2) (p_1^2+ p_2^2) + 2 (R_1 P_1+ R_2 p_2)^2]    dx \nonumber \\
    &&~~~~~\geq  \bar C_1 (\| p_1 \|_2^2+ \| p_2 \|_2^2 )-A_2 (\| w_1 \|_{H^1}^3+\| w_2 \|_{H^1}^3)
    -A_3 (\| w_1 \|_{H^1}^4+\| w_2 \|_{H^1}^4) \nonumber\\
    &&~~~~~\geq - A_2 (\| w_1 \|_{H^1}^3+\| w_2 \|_{H^1}^3) -A_3 (\| w_1 \|_{H^1}^4+\| w_2 \|_{H^1}^4),
    \label{eq10}
 \eea
 where $- \gamma  \int\limits_\rr [(R_1^2+R_2^2) (p_1^2+ p_2^2) + 2 (R_1 P_1+ R_2 p_2)^2]    dx
 \geq -6 \gamma  E^2 (\|p_1 \|_2^2+ \|p_2 \|_2^2) $ is used, and  $~\bar C_1= (C_3-6 \gamma m E^2)/(m+1)$,
 $~ A_2=\max(E_1,E_2)/ (m+1)$ and
 $ ~A_3=\max(C_8,C_9)/(m+1)$  are positive constants. Recalling  that
 $ \langle L_1 p_1,p_1 \rangle+  \langle L_2 p_2,p_2\rangle
      +  2\langle L_3  p_1,p_2\rangle =I+m( \|p_1 \|+ \|p_2 \| )/(m+1)$ we  obtain (\ref{L+})
 where $~ A_1= m \min(1,\Omega)/(1+m)$. This completes the proof of the lemma.
$\hfill\square$

 \noindent
 Finally, the integral term in (\ref{lvar}) can be estimated as
 \bea
 &&\left |-\gamma \int\limits_\rr \left[ \frac{1}{2}(p_1^2+q_1^2+ p_2^2+q_2^2)^2
           + 2 (p_1^2+q_1^2+p_2^2+q_2^2)(p_1 R_1+p_2 R_2) \right]dx \right | \nonumber\\
  &&~~~~~~~~~~~~~~~~
    \leq   \bar D_1 \|w_1 \|_{H^1}( \|w_1 \|_{H^1}^2+\|w_2 \|_{H^1}^2) \nonumber\\
  &&~~~~~~~~~~~~~~~~
           ~~+ \bar D_2 \|w_2 \|_{H^1} ( \|w_1 \|_{H^1}^2+\|w_2 \|_{H^1}^2)
           +\gamma \|w_1 \|_4^4 +\gamma \|w_2 \|_4^4,\nonumber \\
 &&~~~~~~~~~~~~~~~~
    \leq   D_1 \|w_1 \|_{H^1}^3 +D_2 \|w_2 \|_{H^1}^3
           +D_3\|w_1 \|_{H^1}^4 +D_4 \|w_2 \|_{H^1}^4, \label{int.term}
 \eea
 where continuous embedding of $H^1(\rr)$ in $L^4(\rr)$ and in $L^\infty(\rr)$  and
 Young's inequality
 $~ab\leq a^p/p+ b^q/q~\mbox{ with} ~p=3$ and
 $~q=3/2$,  are used, and  $D_i~~(i=1,2,3,4)$ are positive constants.

 \noindent
 \textbf{Proof of Theorem 1.}
 Combining the inequalities (\ref{qiqi}), (\ref{L+}) and (\ref{int.term}),  an upper bound for
 $\Delta L$  is given in terms of $H^1$ norms of the increment
 functions $w_i$ as follows
 \be
 \Delta L(t)\geq g(\| w_1 \|_{H^1})+g(\| w_2 \|_{H^1}), \label{lbound}
 \ee
 where $g(x)=a_1 x^2-a_2 x^3-a_3 x^4$ with positive constants
 \bd
 a_1=\min(C_1,C_2, A_1),~~~ a_2=A_2+\max(D_1,D_2),~~~a_3=A_3+\max(D_3,D_4).
 \ed
 Because $g(0)=0$ and $g(x)\approx a_1 x^2$ near $x=0$, there exists a positive number
 $\epsilon$, $0<\epsilon<\epsilon_0$, such that $g(x)$ increases on $[0,\epsilon_0]$.
 For such an $\epsilon$, the inequalities
 \bd
    \|w_1(0)\|_{H^1}=\| \phi_0(.)-\Phi(.) \|_{H^1} \leq \delta,\;\;\;
    \|w_2(0)\|_{H^1}=\| \psi_0(.)-\Psi(.) \|_{H^1} \leq \delta,
 \ed
 imply that
 \bd
    \Delta L(0)<g(\epsilon)+g(\epsilon)
 \ed
 for sufficiently small $\delta$. As $L(t)$ is invariant with  time, i.e.
 $\Delta L(t)=\Delta L(0)$; from (\ref{lbound}), we have
 \bd
    g(\|w_1(t)\|_{H^1})+g(\|w_2(t) \|_{H^1}) \leq
        \Delta L(t)=\Delta L(0)<g(\epsilon)+g(\epsilon).
 \ed
 By continuity of the function $g$, there is at least a number
 $\epsilon\leq\epsilon_1\leq\epsilon_0$ such  that
 \bd
 \|w_1(t)\|_{H^1}\leq \epsilon_1\leq c_1 \epsilon~~\mbox{and}~~
 \|w_2(t)\|_{H^1}\leq \epsilon_1\leq c_2 \epsilon,
 \ed
 where $t\in [0,\infty)$ and   $c_i~(i=1,2)$ are positive constants.

 Finally for the increment $\eta(x,t)$,  we have to prove that
 $~\displaystyle \|\eta(t)\|_2 \leq c\epsilon$  using the results obtained for
 $\|w_1(t)\|_{H^1}$ and $\|w_2(t)\|_{H^1}$.   In  (\ref{lbound}) we have shown that
 \beas
 \Delta L(t)&=& K+ \frac{c}{2}\int\limits_\rr  \left [ \eta+\frac{2 \beta}{c} (p_1 R_1+ p_2R_2)
        +\frac{\beta}{c}(p_1^2+q_1^2+p_2^2+q_2^2) \right ]^2dx \\
        & \geq & g(\|w_1(t)\|_{H^1})+ g(\|w_2(t)\|_{H^1}) \\
        &&~ +  \frac{c}{2}\int\limits_\rr
        \left [ \eta+\frac{2 \beta}{c} (p_1 R_1+ p_2R_2)
        +\frac{\beta}{c}(p_1^2+q_1^2+p_2^2+q_2^2) \right ]^2dx,
 \eeas
 where
 \beas
 && K=\langle L_0 q_1,q_1 \rangle+\langle L_0 q_2,q_2 \rangle
         +\langle L_1 p_1,p_1 \rangle + \langle L_2 p_2,p_2 \rangle
        + 2 \langle L_3 p_1,p_2 \rangle \\
 &&~~ -\gamma \int\limits_\rr \left [ \frac{1}{2}(p_1^2+q_1^2+ p_2^2+q_2^2)^2
        + 2 (p_1^2+q_1^2+p_2^2+q_2^2)(p_1 R_1+p_2 R_2) \right ]dx.
 \eeas
 For a given $\epsilon>  0$ with $0<\epsilon<\epsilon_0$, the function $g$ is increasing and
 $~g(\|w_i(t)\|_{H^1})>0$ for $~\|w_i(t)\|_{H^1}<c_i \epsilon~~(i=1,2)$. This shows that $K>0$.
 By the invariance property of the functional $L$, $\Delta L(t)=\Delta L(0)$, we have
 \bd
 \int\limits_\rr  \left [ \eta+\frac{2 \beta}{c} (p_1 R_1+ p_2R_2)
        +\frac{\beta}{c}(p_1^2+q_1^2+p_2^2+q_2^2) \right ]^2dx  \leq \frac{4}{c} g(\epsilon).
 \ed
 Using  the inequalities $\displaystyle (a+b)^2\geq \frac{a^2}{2}-b^2$ and
 $(a+b)^2\leq 2(a^2+b^2)$, we find
 \be
    \|\eta(t)\|_2^2\leq \frac{8}{c}g(\epsilon)
    +c_3(\|w_1(t)\|_{H^1}^2+ \|w_2(t)\|_{H^1}^2)+c_4(\|w_1(t)\|_{H^1}^4+ \|w_2(t)\|_{H^1}^2),
 \ee
 where the embedding of $H^1(\rr)$ into $L^2(\rr)$ and $L^4(\rr)$  is used, and $c_3$ and $c_4$ are positive
 constants.  For  some $c>0$, we have $~\displaystyle \|\eta(t)\|_2 \leq c\epsilon$.
 Thus we have proved that solitary waves $(\phi_s, \psi_s, u_s)$  (\ref{solution}) are orbitally stable
 with respect to the small perturbations preserving the $L^2$ norms.

 In order to prove stability of solitary waves with respect to general perturbations, we
 consider  a solitary wave solution $(Q_{1\Omega},Q_{2\Omega})$ which satisfy the
 system (\ref{cs1})
 \beas
 && Q''_{1\Omega}-\Omega Q_{1\Omega}+\gamma (Q_{1\Omega}^2
                + Q_{2\Omega}^2) Q_{1\Omega}=0, ~~~\nonumber \\
 && Q''_{2\Omega}-\Omega Q_{2\Omega}+\gamma (Q_{1\Omega}^2
                + Q_{2\Omega}^2) Q_{2\Omega}=0,                \label{cs3}
 \eeas
 where $\|\phi_0\|_2\neq \|Q_{1\Omega}\|_2$ and $\|\psi_0\|_2\neq \|Q_{2\Omega}\|_2$.
 Then, the functions $~ P_i(x)= Q_{i\Omega}(x /\sqrt \Omega )/\sqrt \Omega~~(i=1,2)$,
 satisfy
 \beas
 && P''_1- P_1+\gamma (P_1^2 + P_2^2) P_1=0,\\
 && P''_2- P_2+\gamma (P_1^2 + P_2^2) P_2=0, \label{cs4}
 \eeas
 where $\|P_i\|_2=\|Q_{i\Omega}\|_2/\sqrt[4]\Omega~~(i=1,2)$. Thus, for the solution
 $(Q_{1\Omega_0}, Q_{2\Omega_0})$ corresponding to $\Omega_0>0$, we have
 $\|P_i\|_2=\|Q_{i\Omega_0}\|_2/\sqrt[4]\Omega_0$. It is possible to choose $\Omega_0$
 such that $\|\phi_0\|_2= \|Q_{1\Omega_0}\|_2$ and $\|\psi_0\|_2= \|Q_{2\Omega_0}\|_2$.
 In the proof of  stability of solitary waves $(Q_{1\Omega}, Q_{2\Omega})$ relative to general perturbations
 that do not preserve $L^2$ norms,
 assuming the initial data obey the inequalities
 $~\| \phi_0(.)- Q_{1\Omega}(.) e^\frac{i c.}{2} \|_{H^1} \leq \delta~$ and
 $~\| \psi_0(.)- Q_{2\Omega}(.) e^\frac{i c.}{2} \|_{H^1} \leq \delta~$,
 the idea is to apply the preceding stability theory
 for  $(Q_{1\Omega_0}, Q_{2\Omega_0})$ and then to use the triangle inequalities
 \bea
   \|  e^{i\theta_1} \phi(.+x_0,t)- Q_{1\Omega}(.) e^\frac{i c.}{2} \|_{H^1} &\leq &
        \|e^{i\theta_1} \phi(.+x_0,t)- Q_{1\Omega_0}(.) e^\frac{ic .}{2}\|_{H^1} \nonumber \\
 &&      +\|Q_{1\Omega_0}(.) - Q_{1\Omega}(.) \|_{H^1}, \label{es1}\\
  \| e^{i\theta_2} \psi(.+x_0,t)-Q_{2\Omega}(.) e^\frac{ic .}{2} \|_{H^1} &\leq &
        \|e^{i\theta_2} \psi(.+x_0,t)-Q_{2\Omega_0}(.) e^\frac{ic .}{F2}\|_{H^1} \nonumber \\
 &&     +\|Q_{2\Omega_0}(.) - Q_{2\Omega}(.) \|_{H^1}. \label{es2}
 \eea
 The first terms in the right hand side of the inequalities (\ref{es1}) and (\ref{es2}) are bounded
 from above by the orbital stability of the solutions $(Q_{1\Omega_0}, Q_{2\Omega_0})$.
 It remains to determine $\delta$ and  to show that $\|Q_{i\Omega_0} - Q_{i\Omega} \|_{H^1}~(i=1,2)$ are also
 small. From the definitions of $Q_{i\Omega}$ and $Q_{i\Omega_0}$ we have
 \bea
  \|Q_{i\Omega} - Q_{i\Omega_0} \|_{H^1}^2&=& \sqrt \Omega
        \int\limits_\rr| ~P_i(x)
                       -\sqrt{\frac{\Omega_0}{\Omega}} ~P_i(\sqrt{\frac{\Omega_0}{\Omega}} ~x)  |^2dx \nonumber\\
 && \!\!\!\!\!\!\!\!+\sqrt{\Omega^3} \int\limits_\rr| ~P_i'(x)
                       -\frac{\Omega_0}{\Omega} ~P_i'(\sqrt{\frac{\Omega_0}{\Omega}} ~x)  |^2dx ~~(i=1,2)\label{qh}.
 \eea
 Using the inequality $(a-\epsilon b)^2\leq 2 \epsilon^2(a-b)^2+2(1- \epsilon)^2 a^2$,
 (\ref{qh}) is rewritten as
 \bea
&&  \|Q_{i\Omega} - Q_{i\Omega_0} \|_{H^1}^2= \sqrt 2 \Omega \left (
        \frac{\Omega_0}{\Omega} \int\limits_\rr| ~P_i(x)
        - ~P_i(\sqrt{\frac{\Omega_0}{\Omega}} ~x)  |^2dx\right.\nonumber\\
  &&~~~~~~~~~~~~~~~~~~~~~~~
   \left.  +( \frac{\Omega_0}{\Omega} -1)^2 \int\limits_\rr  ~P_i^2(x)dx \right )\nonumber\\
 &&~~~~~~~~~~~~~~~~~~~~~~~
 +2\sqrt{\Omega^3} \left (  \frac{\Omega_0^2}{\Omega^2} \int\limits_\rr| ~P_i'(x)
        - P_i'(\sqrt{\frac{\Omega_0}{\Omega}} ~x)  |^2dx\right.\nonumber\\
  &&~~~~~~~~~~~~~~~~~~~~~~~
   \left. +( \frac{\Omega_0}{\Omega} -1)^2 \int\limits_\rr  ~(P_i'(x))^2dx \right ). \label{qh1}
 \eea
 Following the results of Angulo \emph{et. al.} \cite{pava2}, obtained in the study
 of the stability of  solitary waves in the critical case for a generalized Korteweg-de
 Vries  equation and a generalized NLS equation, an upper bound for (\ref{qh1})
 can be given as follows
 \beas
  \|Q_{i\Omega_0} - Q_{i\Omega} \|_{H^1}^2&\leq & G_i~
        (\sqrt[4]{\Omega_0}-\sqrt[4]{\Omega})^2
    +H_i~(\sqrt{\Omega_0} -\sqrt{\Omega})^2~~(i=1,2),
 \eeas
 where
 the fundamental theorem of calculus and Minkowski's inequality are used, and
 the positive constants $G_i$ and $H_i~~(i=1,2)$ are given as
 \bd
     G_i=8\sqrt\frac{\Omega_0}{\Omega}   \left (\|x P_i'\|_2^2
       + \Omega_0 \|x P_i''\|_2^2 \right ),~~~
     H_i= \frac{2}{\sqrt \Omega}   \left (\|P_i\|_2^2
       +  (\sqrt \Omega_0+\sqrt \Omega )^2   \|P_i'\|_2^2 \right ).
 \ed
 We now show that there exists a positive constant $C=C(\Omega_0, P_i)$ such that
 $| \sqrt{\Omega_0}-\sqrt{\Omega}|\leq C \delta$ at least for small values of $\delta$.
 Using the results
 \bd
 \sqrt{\Omega_0}=\frac{\|Q_{i\Omega_0}\|_2^2}{\|P_i\|_2^2}
                =\frac{\|\phi_0\|_2^2}{\|P_1\|_2^2} =\frac{\|\psi_0\|_2^2}{\|P_2\|_2^2},~~~
 \sqrt{\Omega}=\frac{\|Q_{i\Omega}\|_2^2}{\|P_i\|_2^2},
 \ed
 we have
 \beas
  | \sqrt{\Omega_0}-\sqrt{\Omega}|& \leq & \frac{1}{\|P_1\|_2^2}
        \left | \|\phi_0(.)\|_2^2-\|Q_{1\Omega}(.) e^\frac{i c.}{2}\|_2^2  \right |\\
  & \leq & \frac{1}{\|P_1\|_2^ 2} \left (  \delta \|\phi_0(.)\|_2^2
        +(1+\frac{1}{\delta }) \|\phi_0(.) - Q_{1\Omega}(.) e^\frac{i c.}{2}\|_2^2 \right),
 \eeas
 where the inequality $ \left | \|a\|^2- \|b\|^2\right |\leq \|a-b\|^2 +2 \|a\| \|a-b\|$
 and  Young's inequality are used.  Using
 $ \|\phi_0(.) - Q_{1\Omega}(.) e^\frac{i c.}{2}\|_2^2 \leq \delta^2$
 and $\|\phi_0\|_2^2= \sqrt{ \Omega_0}\|P_1\|_2^2 $, we have
  \bd
  | \sqrt{\Omega_0}-\sqrt{\Omega}| \leq   \frac{1}{\|P_1\|_2^ 2}\left ( \delta \sqrt{\Omega_0} \|P_1\|_2^ 2
    +\delta^2+ \delta \right )\leq C(\Omega_0, P_1) \delta
 \ed
 where $C(\Omega_0, P_i)= \sqrt{\Omega_0} +2/\|P_1\|_2^2$.
 The inequality $| \sqrt{\Omega_0}-\sqrt{\Omega}|\leq C \delta$ implies
 $| \sqrt[4]{\Omega_0}-\sqrt[4]{\Omega}|\leq D \delta$ for some positive constant $D$.
 This completes the proof of Theorem 1.

 \noindent
 \textbf{Acknowledgement.} Authors would like to thank Alp Eden for helpful discussions
 at the beginning of this study.

 \end{document}